\documentclass[12pt]{article}
\usepackage{amssymb,amsmath,latexsym}
\usepackage{amsbsy}
\usepackage{amsfonts}
\usepackage{amsopn}
\usepackage{amsthm}
\usepackage{amsmath,verbatim}

\begin{document}

\begin{doublespace}

\newtheorem{remark}{Remark}
\newtheorem{theorem}{Theorem}
\newtheorem{corollary}{Corollary}
\newtheorem{proposition}{Proposition}
\newtheorem{definition}{Definition}
\newtheorem{lemma}{Lemma}
\newtheorem{rem}{Remark}
\newtheorem{defi}{Definition}

\newcommand{\emm}{\mathcal{M}}
\renewcommand{\P}{\mathbb{P}}
\newcommand{\px}{\Bbb{P}_x}
\newcommand{\eq}{\Bbb{E}^{\Bbb{Q}}}
\newcommand{\eqx}{\Bbb{E}^{\Bbb{Q}}_x}
\newcommand{\n}{\underline{n}}
\newcommand{\pf}{\Bbb{P}^\uparrow}
\newcommand{\pfx}{\Bbb{P}_{x}^\uparrow}
\newcommand{\ee}{\mbox{{\bf e}}/\varepsilon}
\newcommand{\ix}{\underline{X}}
\newcommand{\en}{\mbox{\rm I\hspace{-0.02in}N}}
\newcommand{\R}{\mathbb{R}}
\newcommand{\poxi}{\mbox{\raisebox{-1.1ex}{$\mapsto$}\hspace{-.12in}$\xi$}}
\newcommand{\taup}{\mbox{\raisebox{-1.1ex}{$\mapsto$}\hspace{-.12in}$\tau$}}
\newcommand{\ud}{\mathrm{d}}
\newcommand{\ed}{\stackrel{(d)}{=}}
\newcommand{\eqdef}{\stackrel{\mbox{\tiny$($def$)$}}{=}}
\newcommand{\eqas}{\stackrel{\mbox{\tiny$($a.s.$)$}}{=}}
\def\QED{\hfill\vrule height 1.5ex width 1.4ex depth -.1ex \vskip20pt}
\renewcommand{\theequation}{\thesection.\arabic{equation}}

\newcommand{\sfl}{\nabla_1}
\newcommand{\Sf}{\nabla_1_1}
\newcommand{\F}{\mathcal{F}}
\newcommand{\G}{\mathcal{G}}
\newcommand{\PP}{\mathcal{P}}
\newcommand{\LL}{\mathcal{L}}
\newcommand{\TT}{\mathcal{T}}

\def \up{\undeline p}
\def \op{\overline p}

\newcommand{\bE}{\mathbb{E}}
\newcommand{\bP}{\mathbb{P}}
\newcommand{\ind}{\mathbb{I}}

\renewcommand{\qedsymbol}{$\blacksquare$}

\newtheorem{thm}{Theorem}[section]
\newtheorem{defn}{Definition}[section]
\newtheorem{example}[thm]{Example}
\numberwithin{equation}{section}

\newcommand{\D}{{\rm d}}
\def\ee{\varepsilon}
\def\qed{{\hfill $\Box$ \bigskip}}
\def\MM{{\cal M}}
\def\BB{{\cal B}}
\def\LL{{\cal L}}
\def\FF{{\cal F}}
\def\EE{{\cal E}}
\def\QQ{{\cal Q}}

\def\R{{\mathbb R}}
\def\L{{\bf L}}
\def\E{{\mathbb E}}
\def\F{{\bf F}}
\def\P{{\mathbb P}}
\def\N{{\mathbb N}}
\def\eps{\varepsilon}
\def\wh{\widehat}
\def\pf{\noindent{\bf Proof.} }

\title{\Large \bf An optimal stopping problem for fragmentation processes}
\author{Andreas E. Kyprianou\thanks{Department of Mathematical
Sciences, University of Bath, Claverton Down, Bath, BA2 7AY, U.K.
E-mail: a.kyprianou@bath.ac.uk} \,and Juan Carlos Pardo\thanks{CIMAT, Calle Jalisco s/n,
Col. Valenciana, A. P. 402, C.P. 36000, Guanajuato, Gto. MEXICO.
Email: jcpardo@cimat.mx}}
\date{\today}

\maketitle

\begin{abstract}
In this article we consider a toy example of an optimal stopping problem driven by fragmentation processes. 
We show that one can work with the concept of stopping lines to formulate the notion of an optimal stopping problem and moreover, to reduce it to a classical optimal stopping problem for a generalized Ornstein-Uhlenbeck process associated with Bertoin's tagged fragment. We go on to solve the latter using a classical verification technique thanks to the application of aspects of the modern theory of integrated exponential L\'evy processes.
\end{abstract}



\section{Homogenous fragmentation processes}
Fragmentation processes have been the subject of an increasing body of literature and the culmination of this activity has recently been summarised in the recent book of Bertoin \cite{bertfragbook}. Some of the  mathematical roots of fragmentation processes lay with older families of spatial branching processes that have also seen periods of extensive interest such as branching random walks and Crump-Mode-Jagers processes. Irrespective of modern or classical perspectives, such models exemplify the phenomena of random splitting according to systematic  rules and, as stochastic processes, they may be seen as modelling  the growth of special types of multi-particle systems.  The aim of this paper is not to shed light on the, already well understood, intrinsic probabilistic structure of fragmentation processes. Instead we wish to address a previously never treated issue of how a homogenous fragmentation process may be used to drive optimal stopping problems. We give here a toy example and show how the theory of so-called stopping lines allows us to both formulate and solve it by first converting it to a classical optimal stopping problem for an associated generalised Ornstein-Uhlenbeck process.

In order to give a precise statement of our optimal stopping problem, we shall first devote time to recalling some important properties of homogenous fragmentation processes. 

We are interested in the Markov process  $\mathbf{X}:=\{\mathbf{X}(t) : t\geq 0\}$ where $\mathbf{X}(t) =(X_1(t), X_2(t), \cdots)$ that  takes values in 
\[
\mathcal{S}: = \left\{\mathbf{s}=(s_1, s_2, \cdots) : s_1\geq  s_2 \geq \cdots \geq 0 , \, \sum_{i=1}^\infty s_i = 1\right\},
\]
that is to say, the infinite simplex of decreasing numerical sequences with sum equal to 1. The process $\mathbf{X}$ possesses the fragmentation property, to be understood as follows.  Let $t,u\geq 0$. Given that $\mathbf{X}(t)  = (s_1,s_2,\cdots)$, we have that $\mathbf{X}(t+u)$ has the same law as the variable obtained by ranking in decreasing order the collective elements of the the sequences $\mathbf{X}^{(1)}(u), \mathbf{X}^{(2)}(u),\cdots$ where the latter are independent, random mass partitions with values in $\mathcal{S}$ having the same distribution as $s_1\mathbf{X}(u), s_2\mathbf{X}(u), \cdots$, respectively.

It is known that homogenous fragmentation processes can be characterized by a $\sigma$-finite dislocation measure $\nu$ on $\mathcal{S}$ such that $\nu(\{(1,0,\cdots)\})=0$ and 
\[
\int_{\mathcal{S}}(1-s_1)\nu({\rm d}{\bf s})<\infty.
\]
The general theory of fragmentation also allows one to include the possibility of continuous erosion of mass, however, this feature will be excluded in this article. 
Roughly speaking, the dislocation measure specifies the rate at which blocks split so that a block of mass $x$ dislocates into a mass partition $x{\bf s}$, where $\mathbf{s}\in\mathcal{S}$, at rate $\nu({\rm d}{\bf s})$. To be more precise,  let $\mathcal{P}$ be the space of partitions of the natural numbers. Here a partition of $\mathbb{N}$ is a sequence $\pi=(\pi_1, \pi_2, \cdots)$ of disjoint sets, called blocks, such that $\bigcup_i \pi_i = \mathbb{N}$. The blocks of a partition are enumerated in the increasing order of their least element; that is to say $\min \pi_i\leq \min\pi_j$ when $i\leq j$ (with the convention that $\min \emptyset = \infty$). 
Now consider the measure on  $\mathcal{P}$, 
\begin{equation}
\mu(d\pi) = \int_{\mathcal{S}}\varrho_{\bf s}(d\pi)\nu({\rm d}{\bf s}),
\label{implicit}
\end{equation}
 where $\varrho_{\bf s}$ is the law of Kingman's paint-box based on ${\bf s}$ (cf. Chapter 2 of Bertoin \cite{bertfragbook}). It is known that $\mu$ is  an exchangeable partition measure, meaning that it is invariant under the action of finite permutations on $\mathcal{P}$.
It is also known (cf. Chapter 3 of Bertoin \cite{bertfragbook}) that it is possible to construct a fragmentation process  on the space of partitions $\mathcal{P}$ with the help of  a Poisson point process on $(\mathcal{P}\backslash\mathbb{N})\times \mathbb{N}$, say $\{(\pi(t), k(t)) : t\geq 0\}$,  which has intensity measure $\mu\otimes\sharp$, where $\sharp$ is the counting measure. The aforementioned   $\mathcal{P}$-valued fragmentation process is a  Markov process which we denote by  $\Pi = \{\Pi(t) : t\geq 0\}$ such that at all times $t\geq0$ for which an atom $(\pi(t), k(t))$ occurs in the Poisson point process,
$\Pi(t)$ is obtained from $\Pi(t-)$ by partitioning the $k(t)$-th block  into the sub-blocks $(\Pi_{k(t)}(t-) \cap \pi_j(t) : j=1,2,\cdots)$.
Thanks to the properties of the exchangeable partition measure $\mu$ 
it can be shown that, for each $t\geq 0$, the distribution of $\Pi(t)$ is exchangeable and moreover,  blocks of $\Pi(t)$  have asymptotic frequencies in the sense that for each $i\in\mathbb{N}$, 
\[
|\Pi_i(t)|:=\lim_{n\uparrow\infty} \frac{1}{n}\sharp\{\Pi_i(t)\cap \{1,\cdots, n\} \}
\]
exists almost surely. Further,  these asymptotic frequencies, when ranked in decreasing order, say $|\Pi(t)|^\downarrow = (|\Pi(t)|^\downarrow_1, |\Pi(t)|^\downarrow_2, \cdots)$, form a homogenous mass fragmentation process with dislocation measure $\nu$. 
In this sense the process $\mathbf{X}$ and the process $|\Pi|^\downarrow : = (|\Pi(t)|^\downarrow : t\geq 0)$ have the same law.

For future reference we also note from Proposition 2.8 of Bertoin \cite{bertfragbook}  that the ordered asymptotic frequencies of $\pi$ sampled under $\varrho_{\bf s}$, written $|\pi|^\downarrow$, satisfy $|\pi|^\downarrow={\bf s}$ almost surely  and $|\pi_1|$ is a size-biased sample of ${\bf s}$ almost surely. Here `size-biased sample' means that $\varrho_{\bf s}(|\pi_1| = s_i) =s_i$ for $i=1,2,\cdots$.    



Let us introduce the constant 
\[
\underline{ p}: = \inf\left\{  p\in \mathbb{R} : \int_{\mathcal{S}} \left| 1- \sum_{i=1}^\infty s_i^ {1+p} \right| \nu({\rm d}{\bf s}) <\infty\right\}
\]
which is necessarily in $(-1,0]$.
 It is well known that
\[
\Phi( p): = \int_{\mathcal{S}} \left(  1- \sum_{i=1}^\infty s_i^ {1+p} \right)\nu({\rm d}{\bf s}) 
\]
is strictly increasing and concave for $ p\in(\underline{p},\infty)$. Let us assume the following. 

\medskip

{\bf (A1):} If $\underline{p}=0$ then  
$$\Phi'(0+) =\int_{\mathcal{S}}\left(\sum_{i=1}^\infty s_i \log \left(\frac{1}{s_i}\right)\right)\nu(d{\rm s})<\infty.
$$

The function $\Phi$ has a special meaning in the context of the growth of what is commonly referred to as {\it Bertoin's tagged fragment}.   If one considers the process $\xi =\{\xi_t : t\geq 0\}$, where for all $t\geq 0$,
\[
\xi_t := -\log |\Pi_1(t)|,
 \]
  then the underlying Poissonian structure implies that $\xi$ is a subordinator.  Moreover $\Phi$ turns out to be its Laplace exponent meaning that   $$\Phi(p) = -t^{-1}\log \mathbb{E}(e^{-p\xi_t})$$ for all $t,p\geq 0$.  Further, when $\underline{p}<0$, $\xi$ has finite mean, that is to say $\Phi'(0+)<\infty$ and the same is true when $\underline{p}=0$ thanks to (A1). 
Note that the L\'evy measure, $\mu$, associated to $\xi$ is given through the dislocation measure $\nu$ via the formula 
  \[
  m({\rm d}x)  = e^{-x}\sum_{i=1}^\infty\nu(-\log s_i \in {\rm d}x)
  \]
  for $x>0$ and, in particular, it necessarily satisfies $\int_{(0,\infty)}(1\wedge x)m({\rm d}x)<\infty$.

On account of the fact that $\Pi(t)$ is an exchangeable partition, it also follows that, given $|\Pi(t)|^\downarrow$, the distribution of $|\Pi_1(t)|$ is that of a size-biased pick. This observation gives us the so-called {\it Many-to-one Principle} that for all positive, bounded measurable $f$,
\begin{equation}
\mathbb{E}\left[\sum_{n\geq 1} f(|\Pi_n(t)|)|\Pi_n(t)| \right]=  \mathbb{E}[f(e^{-\xi_t})].
\label{M21}
\end{equation}
The above formula, or rather a slightly more elaborate version of it, turns out to be a key ingredient in formulating optimal stopping problems for fragmentation processes.

\section{Stopping lines} \label{S: stop line}

The concept of a {\it stopping line } was introduced by Bertoin \cite{bertfragbook} in the context of fragmentation processes, capturing in its definition the essence of earlier ideas on stopping lines for various classes of spatial branching processes coming from the work of Neveu \cite{neveu},  Jagers \cite{jagers}, Chauvin \cite{chauvin} and Biggins and Kyprianou \cite{BK97}.  Roughly speaking a stopping line plays the an analogous role to that of a stopping time for single particle Markov processes. 
We give below a more concise definition.

Recall that for each integer $i \in \mathbb{N}$ and $s \in \R_+$ we denote by $B_i(s)$ the block of $\Pi(s)$ which contains $i$, with the convention that $B_i(\infty)=\{i\}$. For comparison, recall that  $\Pi_i(t)$ is the $i$th block by order of least element. We write 
\[
 \G_i(t) = \sigma (B_i(s) , s\le t)
\]
for the sigma-field generated by the history of the block containing $i$ up to time $t.$ 

\begin{definition}\rm  
 We call a stopping line a family $\ell=(\ell(i) , i\in\mathbb{N})$ of random variables with values in $[0, \infty]$ such that for each $i \in \mathbb{N}$,
\begin{itemize}
 \item[(i)] $\ell(i)$ is a $(\G_i(t))$-stopping time.
 \item[(ii)] $\ell(i)=\ell(j)$ for every $j \in B_i(\ell(i)).$
\end{itemize}
\end{definition}

For instance, first passage times such as $\ell(i) = \inf \{ t\ge 0 : |B_i(t)| \le a \}$ for a fixed level $a\in (0,1)$ define a stopping line.  

\medskip 

The key point is that it can be checked that the collection of blocks $\Pi(\ell)=\{ B_i(\ell(i)) \}_{i \in \mathbb{N}}$ is a partition of $\mathbb{N}$ which we denote by $\Pi(\ell) = (\Pi_1(\ell), \Pi_2(\ell),\cdots),$ where as usual the enumeration is by order of least element.

Observe that because $B_i(\ell(i))=B_j(\ell(j))$ when $j \in B_i(\ell(i))$, the set $\{ B_i(\ell(i)) \}_{i \in \mathbb{N}}$ has repetitions. Hence $(\Pi_1(\ell), \Pi_2(\ell),\cdots)$ is simply a way of enumerating each element only once by order of discovery. In the same way the $\ell(j)$'s can be enumerated as $\ell_i, i=1,\cdots $ such that for each $i, \ell_i$ corresponds to the stopping time of $\Pi_i(\ell).$ Note that $\ell(1) = \ell_1$. For convenience, for each $t\geq 0$, we may write $\{\overleftarrow{\Pi}_i(s): s\leq t\}$ to mean ancestral evolution of the block $\Pi_i(t)$.

The following lemma gives us a  generalisation of (\ref{M21}).



\begin{lemma}[Many-to-one principle]\label{L:many to one}
Let $\ell$ be a stopping line Then, for any non-negative, bounded measurable $f$ and family of random variables $(A_{i} (\ell):i\geq 1)$ such that for each $i$, $A_{i}(\ell) = F(\{\overleftarrow{\Pi}_i(s): s\leq \ell_i\})$, where $F$ is a non-negative measurable function of $\{\overleftarrow{\Pi}_i(s) : s\leq \ell_i\}$. Then
\begin{equation}\label{E:many-to-one}
 \mathbb{E} (\sum_i |\Pi_i(\ell)| f(A_i(\ell), \ell_i)) = \mathbb{E}(f(A_1 (\ell),\ell_1) \mathbf{1}_{\{ \ell_1 <\infty \}}).
\end{equation}
\end{lemma}
{\bf Proof.} The proof is essentially the same as the proof of (\ref{M21}); cf. Lemma 2 of \cite{BHK}.
Indeed, since $\Pi(\ell)$ is an exchangeable partition, the sum on the left hand side can be understood to indicate that the pair $(A_1(\ell),\ell_1)$ is a size-biased pick from the sequence $((A_i(\ell),\ell_i), i\geq 1)$. The indicator function in the right-hand side comes from the possibility that there is some dust in $\Pi(\ell)$ 
\QED

\section{An optimal stopping problem}

The main objective of this paper is to solve the following toy optimal stopping problem for the homogenous fragmentation process described above.
We want to find both the value and an optimizing strategy for the following quantity
\begin{equation}
 V(c) = \sup_\ell \mathbb{E}\left[ 
\sum_{n\geq 1} \left(\int_0^{\ell_n} e^{-\gamma\theta s} |\overleftarrow{\Pi}_n(s)|^{-\gamma} {\rm ds} +c  \right)|\Pi_n(\ell)|^{1+\gamma}e^{-q \ell_n}
\right],
\label{frag-opt-stop}
\end{equation}
where $\gamma,\theta>0$, $q\geq 0$ and the supremum is taken over all stopping lines $\ell$.

To give a meaning to the above optimal stopping problem we can think of the following interpretation. A commodity is crushed and sold. For the sake of argument we may consider this to be rock. Buyers of the rock may order it crushed according to any of the strategies $\ell$ (one may think of this as ordering gravel from the original rock source which has been crushed to a particular specification). The seller of the crushed ensemble $\Pi(\ell)$ prices as follows. For every fragment $\Pi_n(\ell)$ in the crushed ensemble $\Pi(\ell)$, the buyer must pay a premium which is the sum of $c$ monetary units and a cost which depends on the process of how the individual fragment $\Pi_n(\ell)$ was formed; that is $\int_0^{\ell_n} e^{-\gamma\theta} |\overleftarrow{\Pi}_n(s)|^{-\gamma}{\rm d}s$. Note that fragments which are crushed to a very fine degree in a 
short amount of time will attract a large premium. At the same, the seller discounts the premium on individual fragments according to their size through the factor $|\Pi_n(\ell)^{1+\gamma}|$. The smaller a block the larger the discount is offered to offset the large crushing premium. Finally, as the pricing mechanism depends on the evolution of the crushing process through time, a exponential discounting is imposed at rate $q$ to accommodate for the net present value of monetary value.

The first and most important observation in handing (\ref{frag-opt-stop}) is that, with the help of the Many-to-one principle  we are able to convert this optimal stopping problem for fragmentation processes into an optimal stopping problem for a single particle Markov process, which, in this case, turns out to be a generalised Ornstein-Uhlenbeck process. The relevant definition of this process is  as follows. 

Suppose that $Y=\{Y_t :t\geq 0\}$  is a spectrally positive L\'evy process with Laplace exponent given by
\[
 \psi(u) = \log\mathbb{E}(e^{-u Y_1}) = \theta u - \Phi(u)
\]
for $u\geq 0$. That is to say, $Y$ is equal in law to the difference of the subordinator describing Bertoin's tagged fragment and a linear drift with rate $\theta$. (Note that in particular, $Y$ is a process with bounded variation paths). Then, for each $c>0$ define the stochastic process $Z^c = \{Z^c_t : t\geq 0\}$ by
\[
 Z^c_t = e^{-\gamma Y_t} \left( \int_0^t e^{\gamma Y_s}{\rm d}s +c \right), \,\, t\geq 0.
\]

On account of the fact that $Y$ is spectrally positive, it follows that $Z^c$ experiences only negative jumps in its path. It is also easy to see that $Z^c$ is a Markov process. Indeed, writing $Y_{t+s} = Y_t + \widetilde{Y}_s$ where $\widetilde{Y} = \{\widetilde{Y}_u: u\geq 0\}$ is independent of $\{Z^c_u: u\leq t\}$ and has the same law as $Y$, we have
\begin{eqnarray*}
Z^c_{t+s} &=&  e^{-\gamma \widetilde{Y}_s}e^{-\gamma Y_t} \left( \int_0^t e^{\gamma Y_s}{\rm d}s  +c + e^{\gamma Y_t}\int_0^s e^{\gamma  \widetilde{Y}_u}{\rm d}u \right)\\
&=& e^{-\gamma \widetilde{Y}_s} \left( e^{-\gamma Y_t}\left(\int_0^t e^{\gamma Y_s}{\rm d}s  +c\right) + \int_0^s e^{\gamma  \widetilde{Y}_u}{\rm d}u \right)\\
&=& \widetilde{Z}^{Z^c_t}_s
\end{eqnarray*}
where, for each $c>0$, $\widetilde{Z}^c$ is independent of $\{Z^c_u: u\leq t\}$ and has the same law as $Z^c$.

Henceforth let us assume that

\bigskip

{\bf (A2):} $\theta>\Phi'(0+)$.

\bigskip

\noindent in which case it follows that $\mathbb{E}(Y_1) = - \theta + \Phi'(0+)<0 $ and hence $\lim_{t\uparrow\infty}Y_t = -\infty$ almost surely. In that case it is known (cf. \cite{MaZ}) that 
\[
\int_0^\infty e^{\gamma Y_t}{\rm d}t<\infty,
\]
which implies that
 $\lim_{t\uparrow\infty} Z^c_t$ exists and is equal to $+\infty$ almost surely.

\bigskip

The next lemma shows how we may reformulate (\ref{frag-opt-stop}) in terms of the process $Z^c$.

\begin{lemma}\label{many2one-optstop}
 For all $\gamma, \theta, c>0$ and $q\geq 0$,
\begin{equation}
 V(c) = \sup_{\tau}\mathbb{E}[e^{-(q+\theta\gamma)\tau} Z^c_\tau]
 \label{OU-opt-stop}
\end{equation}
where the supremum is taken over the family of stopping times  for the stochastic process $Y$.
\end{lemma}

{\bf Proof.}
A standard argument using the  Many-to-one Principle shows us that  for a given stopping line $\ell$ we have
\begin{eqnarray*}
 \lefteqn{\mathbb{E}\left[ 
\sum_{n\geq 1} \left(\int_0^{\ell_n} e^{-\gamma\theta s} |\overleftarrow{\Pi}_n(s)|^{-\gamma} {\rm ds} +c  \right)|\Pi_n(\ell)|^{1+\gamma}e^{-q \ell_n}
\right]}&& \\
&&= \mathbb{E}\left[ 
 \left(\int_0^{\ell_1} e^{-\gamma\theta} e^{\gamma\xi_s} {\rm ds} +c  \right)e^{-\gamma\xi_{\ell_1}}e^{-q \ell_1}\mathbf{1}_{\{\ell_1<\infty\}}
\right]\\
&&= \mathbb{E}\left[ 
 \left(\int_0^{\ell_1} e^{\gamma Y_s} {\rm ds} +c  \right)e^{-\gamma Y_{\ell_1}}e^{-(q+\theta\gamma) \ell_1}
 \mathbf{1}_{\{\ell_1<\infty\}}
\right]\\
&&= \mathbb{E}\left[ 
Z_{\ell_1}e^{-(q+\theta\gamma) \ell_1}
\mathbf{1}_{\{\ell_1<\infty\}}
\right]
\end{eqnarray*}
where we recall that $\xi=\{\xi_t: t\geq 0\}$ plays the role of the tagged fragment and left hand and right hand sides of this series of equalities are understood to be infinite in value at the same time.
The indicator can be removed in the final equality as soon as we can show that 
\[
\lim_{t\uparrow\infty} e^{-(q+\theta\gamma)t }Z^c_t = 0.
\]
This is easily shown however by writing $\lambda = q+\theta\gamma$ and observing that
\[
e^{-\lambda t}Z^c_t \leq  e^{-\kappa(\lambda)Y_t - \lambda t} e^{-(\gamma-\kappa(\lambda))Y_t} \left(\int_0^\infty e^{\gamma Y_s} {\rm d}s + c\right),
\]
where $\kappa(\lambda)$ the largest solution of the equation $\psi(u)=\lambda$. Note that since $\psi$ is convex, tending to infinity at infinity (cf. Chapter 8 of \cite{kyp book}), $\kappa(\lambda)$ is finite.
Recall that the integral on the right hand side is almost surely convergent. Moreover since $\exp\{-\kappa(\lambda)Y_t -\lambda t\}$ is a positive martingale, and hence almost surely convergent, and since,  by the definition of $\lambda$, $\kappa(\lambda)>\gamma$, it follows that 
 $\lim_{t\uparrow\infty} e^{-\lambda t}Z^c_t =0$ almost surely.
\QED

We remark that a family of optimal stopping problems for generalised Ornstein-Uhlenbeck similar to (\ref{OU-opt-stop}) has been considered by Gapeev \cite{Ga} for  compound Poisson processes  with  exponential jumps and by Ciss\'e et al. \cite{cisse} for spectrally negative L\'evy processes.

\section{Solution of the optimal stopping problem}

It turns out that, like many optimal stopping problems, the optimal strategy in (\ref{frag-opt-stop}) boils down to first passage over a threshold of an auxiliary process. In this case, when considering our optimal stopping problem in the form (\ref{OU-opt-stop}) it will turn out to be optimal to stop when the generalized OU process $Z^c$ crosses an appropriately chosen level $b^*$. 

Define for $b\geq c$,
\[
\tau^c_{(b,\infty)} = \inf\{t>0 : Z^c_t >b\}.
\]
Note in particular that, thanks to the fact that $Z^c$ is skip-free upwards, we have that $Z^c_{\tau_{[0,\infty)}} = b$ on the event $\{\tau_{(b,\infty)}<\infty\}$ which, in turn, occurs with probability 1 thanks to the assumption (A2). Hence,  for any $\lambda>0$,
\[
\mathbb{E}(e^{-\lambda \tau^c_{(b,\infty)}  } Z^c_{ \tau^c_{(b,\infty)} }) = 
b\mathbb{E}(e^{-\lambda \tau^c_{(b,\infty)}   }).
\]

Recall that, for $\lambda\geq 0$, $\kappa(\lambda)$ is the largest solution of the equation $\psi(u)=\lambda$. 
Define, for all $t\geq 0$, $\mathcal{F}_t = \sigma(Y_s: s\leq t)$ and consider the exponential change of measure
\begin{equation}
\frac{\ud \mathbb{P}^{\kappa(\lambda)}}{\ud \mathbb{P}}\bigg|_{\mathcal{F}_t}=e^{-\kappa(\lambda) Y_t-\lambda t},\qquad \textrm{for } \lambda\ge 0.
\label{esscher}
\end{equation}
Under $\mathbb{P}^{\kappa(\lambda)}$, the process $Y$ is still a spectrally positive and its Laplace exponent, $\psi_{\kappa(\lambda)}$ satisfies the relation
\[
\psi_{\kappa(\lambda)}(u)=\psi(\kappa(\lambda)+u)-\lambda, \qquad\textrm{ for}\quad u\ge 0. 
\]
See for example Chapter 8 of \cite{kyp book}
 for further details on the above remarks. 
 Note in particular that it is easy to verify that $\psi_{\kappa(\lambda)}'(0+)>0$ and hence the process $Y$ under $\mathbb{P}^{\kappa(\lambda)}$ drifts to $-\infty$. According to earlier discussion, this guarantees that also under $\mathbb{P}^{\kappa(\lambda)}$, the process $Z^c$ also drifts to $+\infty$.

\begin{lemma}\label{first passage LT}
Suppose that $\lambda\geq 0$ and that $\kappa(\lambda)>\gamma$, then for all $b\geq c$,
\[
\mathbb{E}(e^{-\lambda \tau^c_{(b,\infty)}}) = \frac{\mathbb{E}^{\kappa(\lambda)} [ (c + I_\infty)^{\frac{\kappa(\lambda)}{\gamma}}]}{\mathbb{E}^{\kappa(\lambda)} [ (b + I_\infty)^{\frac{\kappa(\lambda)}{\gamma}}]},
\]
where 
\[
I_\infty = \int_0^\infty e^{\gamma Y_s}{\rm d}s.
\]
\end{lemma}

{\bf Proof.} The proof relies on the Lamperti representation of  the positive self-similar Markov process associated to $Y$ and follows similar arguments as in the proof of Theorem 1 in Enriquez et al. \cite{ESY}. Recall, from the Lamperti representation (cf. \cite{La}), that the process $X=\{X_t,t\ge 0\}$  defined as follows
\[
X_t=\exp\Big\{Y_{I^{-1}(t)}\Big\},\quad \textrm{ where }\quad I^{-1}(t)=\inf\left\{s\ge 0:\,I_s>t\right\} \textrm{ and } I_s=\int_0^s e^{\gamma Y_u} \ud u,
\]
is a positive self-similar Markov process starting from $1$  with self-similar index equal to $\gamma>0.$ Let
\[
\sigma(b,c)=\inf\{t\ge 0 : \big(X_t\big)^{\gamma}\le \left(t +c\right)/b\},
\]
and note from the Lamperti representation,  that $\sigma(b,c)=I_{\tau^{c}_{(b,\infty)}}$.

For any almost surely stopping time $T$ with respect to $\{\mathcal{F}_t : t\geq 0\}$, we have
\begin{equation}\label{equ1}
\begin{split}
I_\infty&=\int_0^\infty e^{\gamma Y_u} \ud u=\int_0^T e^{\gamma Y_u} \ud u+e^{\gamma Y_T}\int_0^\infty e^{\gamma Y'_u} \ud u,\\
&=I_T+e^{\gamma Y_T}I_\infty^\prime,
\end{split}
\end{equation}
where $Y^\prime_t:=Y_{t+T}-Y_{T}$, for $t\ge 0$ and $I'_\infty$ is defined by the second equality. From the strong Markov property of $Y$, we observe that the random variable $I_\infty^\prime$ is independent of  $ (I_T,e^{\gamma Y_T})$ and  has the same law as $I_\infty$. 

Next note that $I_{\tau^c_{(b,\infty)}} = \sigma(b,c)$. Indeed, as $Z^c$ has no positive jumps we have  
\begin{equation}
e^{\gamma Y_{\tau^c_{(b,\infty)}}} = \frac{1}{b}(I_{\tau^c_{(b,\infty)}} +c),
\label{eq1}
\end{equation}
from which it follows that $I_{\tau^c_{(b,\infty)}} \geq \sigma(b,c)$. Moreover, since $X_{\sigma(b,c)}^\gamma = (\sigma(b,c) +c)/b$, that is to say,
\[
b= e^{-\gamma Y_{I^{-1}(\sigma(b,c))}} (I_{I^{-1}(\sigma(b,c))} + c),
\]
we also have that $I^{-1}(\sigma(b,c))\geq \tau^c_{(b,\infty)}$, or equivalently $\sigma(b,c)\geq I_{\tau^c_{(b,\infty)}}$.

Using  (\ref{equ1}) with $T = \tau^c_{(b,\infty)}$, we now have 
\[
\begin{split}
c+I_\infty&=c+I_{\tau^{c}_{(b,\infty)}}+e^{\gamma Y_{\tau^{c}_{(b,\infty)}}}I_\infty^\prime\\
&=c+\sigma(b,c)+\frac{1}{b}\left(\sigma(b,c) +c\right)I_\infty^\prime\\
&=(c+\sigma(b,c))\left(1+\frac{1}{b}I_\infty^\prime\right),
\end{split}
\]
which implies 
\begin{equation}
\mathbb{E}^{\kappa(\lambda)}[(1+b^{-1}I_\infty)^{s}]\mathbb{E}^{\kappa(\lambda)}[(c+\sigma(b,c))^{s}]=\mathbb{E}^{\kappa(\lambda)}[(c+I_\infty)^{s}] \qquad \textrm{for any} \quad s\in \R,
\label{gives the result}
\end{equation}
where we understand both sides of the equality to be infinite simultaneously.
In order to complete the proof, we need to show that for   $s = \frac{\kappa(\lambda)}{\gamma}$ and $a\geq 0$, the  quantity
$\mathbb{E}^{\kappa(\lambda)}[(a+I_\infty)^{s}]$ is finite. In that case, it follows from (\ref{esscher}), (\ref{eq1})  and (\ref{gives the result}) that for $\lambda\geq 0$ 
\begin{eqnarray}
\mathbb{E}[e^{-\lambda \tau^c_{(b,\infty)}}] &=& \mathbb{E}^{\kappa(\lambda)} (e^{\kappa(\lambda) Y_{\tau^c_{(b,\infty)}}}) \notag\\
& =& b^{-\frac{\kappa(\lambda)}{\gamma}}\mathbb{E}^{\kappa(\lambda)}[(c+\sigma(b,c))^{\frac{\kappa(\lambda)}{\gamma}}] \notag\\
& =& \frac{\mathbb{E}^{\kappa(\lambda)}[(c+I_\infty)^{\frac{\kappa(\lambda)}{\gamma}}]}{\mathbb{E}^{\kappa(\lambda)}[(b+I_\infty)^{\frac{\kappa(\lambda)}{\gamma}}]}\label{similar one later}
\end{eqnarray}
as required.

Since, for $s\geq 1$
\[
\mathbb{E}^{\kappa(\lambda)}[(a+I_\infty)^s]\le 2^{s -1}(a^s+\mathbb{E}^{\kappa(\lambda)}[I_\infty^s]),
\]
it suffices to investigate the finiteness of $\mathbb{E}^{\kappa(\lambda)}[I_\infty^s]$.
According to  Lemma 2.1 in Maulik and Zwart \cite{MaZ} the expectation $\mathbb{E}^{\kappa(\lambda)}[I_\infty^{s}]$ is finite for all $s\ge 0$ such that $-\psi_{\kappa(\lambda)}(-\gamma s)>0$. Since $\psi_{\kappa(\lambda)}(-\gamma s)$ is well defined for $\kappa(\lambda)-\gamma s\ge 0$, 
then a straightforward computation gives us that $\mathbb{E}^{\kappa(\lambda)}[I_\infty^s]<\infty$ for $s\in [0,\frac{\kappa(\lambda)}{\gamma}]$. 
\QED

The next lemma concerns how we can choose the optimal threshold, $b^*$.

\begin{lemma}Suppose that $\kappa(\lambda)>\gamma$ and define  the function 
\[
f(b) = \frac{1}{b}\frac{\mathbb{E}^{\kappa(\lambda)} [ (b + I_\infty)^{\frac{\kappa(\lambda)}{\gamma}}]}{\mathbb{E}^{\kappa(\lambda)} [ (b + I_\infty)^{\frac{\kappa(\lambda)}{\gamma} - 1}]}.
\]
Then $f$ is continuous, strictly monotone decreasing, satisfies $f(0+) = \infty$ and $f(\infty)  =1$ and there exists a unique solution, denoted by $b^*$, to the functional equation
\[
f(b) = \frac{\kappa(\lambda)}{\gamma}.
\]
\end{lemma}
{\bf Proof.} It is clear from the definition of $f$ that $f(0+)=\infty$ and $f(\infty) =1$.
Let $b<a$ and note that monotonicity of $f$ follows on account of the fact that
\[
\begin{split}
f(a)&=\frac{1}{a}\frac{\mathbb{E}^{\kappa(\lambda)}[(a+I_\infty )^{\frac{\kappa(\lambda)}{\gamma}}]}{\mathbb{E}^{\kappa(\lambda)}[(a+I_\infty )^{\frac{\kappa(\lambda)}{\gamma}-1}]}\frac{\mathbb{E}^{\kappa(\lambda)}[(b+I_\infty )^{\frac{\kappa(\lambda)}{\gamma}}]}{\mathbb{E}^{\kappa(\lambda)}[(b+I_\infty )^{\frac{\kappa(\lambda)}{\gamma}}]}\\
&=\frac{1}{a}\frac{1}{\mathbb{E}^{\kappa(\lambda)}[e^{-\lambda \tau^b_{(a,\infty)}}]}\frac{\mathbb{E}^{\kappa(\lambda)}[(b+I_\infty )^{\frac{\kappa(\lambda)}{\gamma}}]}{\mathbb{E}^{\kappa(\lambda)}[(a+I_\infty )^{\frac{\kappa(\lambda)}{\gamma}-1}]}\\
&=\frac{b}{a}\frac{f(b)}{\mathbb{E}^{\kappa(\lambda)}[e^{-\lambda \tau^b_{(a,\infty)}}]}\frac{\mathbb{E}^{\kappa(\lambda)}[(b+I_\infty )^{\frac{\kappa(\lambda)}{\gamma}-1}]}{\mathbb{E}^{\kappa(\lambda)}[(a+I_\infty )^{\frac{\kappa(\lambda)}{\gamma}-1}]}\\
&=\frac{b}{a}\frac{f(b)}{\mathbb{E}^{\kappa(\lambda)}[e^{-\lambda \tau^b_{(a,\infty)}}]}\mathbb{E}^{\kappa(\lambda)}[e^{-\lambda \tau^b_{(a,\infty)}-\gamma Y_{\tau^b_{(a,\infty)}}}]\\
&\le f(b),
\end{split}
\] 
where the penultimate equality is the result of a computation similar to (\ref{similar one later}).
\QED

Define  the pair $(\tau^*, V^*)$ such that $\tau^*= \tau^c_{(b^*,\infty)}$ and 
$$
V^*(c)  = \mathbb{E}(e^{-\lambda \tau^*}   Z^c_{ \tau^* }).
$$
Our main result giving the solution to our optimal stopping problem can now be stated as follows.

\begin{theorem}\label{main} Suppose that  $q,\gamma>0$ and recall that (A1) and (A2) are in force.
The value function $V$ in (\ref{OU-opt-stop}) is equal to $V^*$ with $\lambda = q+\theta\gamma$ and this value is obtained by the optimal strategy $\tau^*$. In particular this implies that an  optimal stopping line strategy is given by $(\ell^*(n):n\in\mathbb{N})$ where
\[
\ell^*(n) =  \inf\left\{ 
t>0 : 
\left(\int_0^{t} e^{-\gamma\theta s} |B_n(s)|^{-\gamma} {\rm ds} +c  \right)|B_n(t)|^{\gamma} > b^*\right\}.
\]
\end{theorem}

\section{Proof of Theorem \ref{main}}
Note that, by exchangeability,  $\ell^*(1)$ characterises $\ell^*(n)$ for all $n\in\mathbb{N}$. At the same time we have $\ell^*(1) = \ell^*_1$ and hence the proof is complete as soon as we show that 
\[
\ell^*_1 = \inf\left\{ 
t>0 : 
\left(\int_0^{t} e^{-\gamma\theta s} |B_1(s)|^{-\gamma} {\rm ds} +c  \right)|B_1(t)|^{\gamma} > b^*\right\},
\]
which is equivalent to showing the proof of (\ref{OU-opt-stop}) is 
given by the pair $(\tau^*,V^*)$.

Following a classical verification technique, it suffices to show that \begin{itemize}
\item[(i)] the limit $\lim_{t\uparrow\infty} e^{-\lambda t}Z^c_t$ exists almost surely and is finite, 
\item[(ii)] $V^*(c)\geq c$ for all $c> 0$ and
\item[(iii)]  for all $c>0$, $\{e^{-\lambda t} V^*(Z^c_{t}): t\geq 0\}$ is a  right continuous super-martingale.
\end{itemize}
See for example Theorem 9.1 in Kyprianou \cite{kyp book}.

The requirement  (i) was demonstrated in the proof of Lemma \ref{many2one-optstop}.  
In order to show (ii) and (iii) we need to introduce  the following function, 
\[
\widetilde{V}(c) =  b^*\frac{\mathbb{E}^{\kappa(\lambda)} [ (c + I_\infty)^{\frac{\kappa(\lambda)}{\gamma}}]}{\mathbb{E}^{\kappa(\lambda)} [ (b^* + I_\infty)^{\frac{\kappa(\lambda)}{\gamma}}]},
\]
defined for all $c>0$.
Recall that, since  $\lambda = q + \gamma\theta$, $\kappa(\lambda)>\gamma$ and hence with the help of Lemma \ref{first passage LT}, $V^*(c) = \widetilde{V}(c)$ for all $c\leq b^*$ and otherwise $V^*(c) = c$ for $c> b^*$. There is also continuity at $b^*$ in $V^*$ thanks to regularity of $(b^*,\infty)$  for $Z^{b^*}$, i.e. $\mathbb{P}(\tau^{b^*}_{(b^*,\infty)} =0) = 1$.  
Moreover, on account of the fact that $\frac{\kappa(\lambda)}{\gamma}>1$, we easily see that $\widetilde{V}$ is a convex, continuously differentiable function. 
In particular, thanks to the definition of $b^*$, it is also clear that  $\widetilde{V}'(b^*) = 1$ thus making $c$ a tangent line to $\widetilde{V}(c)$, touching it at $b^*$, and hence $\widetilde{V}(c)\geq c$ for all $c>0$. Specifically we note that there is both {\it continuous and smooth pasting} at $b^*$. It follows that $V^*$ is convex,  continuously differentiable and $V^*(c)\geq c$. The latter is the requirement (ii).

 It thus remains to show (iii). We do this through the following three lemmas.

\begin{lemma}
For all $c>0$, $\{e^{-\lambda t}\widetilde{V}(Z^c_t): t\geq 0\}$ is a martingale. 
\end{lemma}
{\bf Proof.}  Note that 
\begin{eqnarray}
\lefteqn{\mathbb{E}(e^{-\lambda t} \widetilde{V}(Z^c_t))}&&\nonumber\\
&=& b^*\mathbb{E}\left(e^{-\lambda t}\frac{\mathbb{E}^{\kappa(\lambda)} [ (Z^c_t + \widetilde{I}_\infty)^{\frac{\kappa(\lambda)}{\gamma}} |\mathcal{F}_t]}{\mathbb{E}^{\kappa(\lambda)} [ (b^* + I_\infty)^{\frac{\kappa(\lambda)}{\gamma}}]}\right)\nonumber\\
&=&b^*\mathbb{E}\left(e^{-\lambda t} e^{-\kappa(\lambda) Y_t}\frac{\mathbb{E}^{\kappa(\lambda)} [ (c + \int_0^t e^{\gamma Y_s}{\rm d}s+ e^{\gamma Y_t}\widetilde{I}_\infty)^{\frac{\kappa(\lambda)}{\gamma}} | \mathcal{F}_t]}{\mathbb{E}^{\kappa(\lambda)} [ (b^* + I_\infty)^{\frac{\kappa(\lambda)}{\gamma}}]}\right)\nonumber\\
&=&b^*\mathbb{E}^{\kappa(\lambda)}\left(\frac{\mathbb{E}^{\kappa(\lambda)} [ (c + \int_0^t e^{\gamma Y_s}{\rm d}s+ e^{\gamma Y_t}\widetilde{I}_\infty)^{\frac{\kappa(\lambda)}{\gamma}} | \mathcal{F}_t]}{\mathbb{E}^{\kappa(\lambda)} [ (b^* + I_\infty)^{\frac{\kappa(\lambda)}{\gamma}}]}\right)\nonumber\\
&=&b^*\frac{\mathbb{E}^{\kappa(\lambda)} [ (c + I_\infty)^{\frac{\kappa(\lambda)}{\gamma}}]}{\mathbb{E}^{\kappa(\lambda)} [ (b^* + I_\infty)^{\frac{\kappa(\lambda)}{\gamma}}]}\nonumber\\
&=& \widetilde{V}(c),\label{constantmean}
\end{eqnarray}
where $\widetilde{I}_\infty$ is independent of $\mathcal{F}_t$ and has the same distribution as $I_\infty$.
Now note by the Markov property that 
\begin{equation}
\mathbb{E}(e^{-\lambda(t+s)} \widetilde{V}(Z^c_{t+s})|\mathcal{F}_t)
 = e^{-\lambda t}\mathbb{E}(e^{-\lambda s} \widetilde{V}(\widetilde{Z}^{Z^c_t}_s)|\mathcal{F}_t) = e^{-\lambda t}\widetilde{V}(Z^c_t)
 \label{a similar computation}
\end{equation}
where, given $c>0$, $\{\widetilde{Z}^c_t: t\geq 0\}$ is a copy of $Z^c$ which is independent of $\mathcal{F}_t$.
 \QED

Next,  note that it is straightforward to show that the process $\{Z^c_t:t\geq 0\}$ is described by the stochastic differential equation
\[
Z^c_t = c+  \int_0^t (1 + \gamma\theta Z^c_{s-}){\rm d} s -\sum_{0<s\leq t}Z^c_{s-}(1-e^{-\gamma\Delta \xi_s}),
\]
where we recall that $\xi = \{\xi_t : t\geq 0\}$ is Bertoin's tagged fragment and, for all $t>0$, $\Delta \xi_t = \xi_t - \xi_{t-}$. Using standard calculus for bounded variation, right-continuous stochastic processes, one deduces that for for any  $C^1(0,\infty)$ function $f$, 
\begin{equation}
e^{-\lambda t} f(Z^c_{t}) - f(c) = \int_0^t e^{-\lambda s} (\mathcal{L}-\lambda)f(Z^c_{s-}){\rm d}s + M^f_t,
\label{itodecomp}
\end{equation}
where $c>0$, $t\geq 0$, 
\[
\mathcal{L}f(x)  = (1+\gamma\theta x)f'(x) 
+ \int_{(0,\infty)} \left\{f(e^{-\gamma y}x)-f(x)\right\}m({\rm d}y),
\]
and the process $\{M^f_t: t\geq 0\}$ is a martingale given by the compensated process
\[
M^f_t = \sum_{s\leq t} \left\{f(e^{-\gamma \Delta\xi_t} Z^c_{s-} ) - f(Z^c_{s-})\right\}
- \int_0^t\int_{(0,\infty)} \left\{f(e^{-\gamma y} Z^c_{s-} ) - f(Z^c_{s-})\right\}m({\rm d}y){\rm d}s
\]
for $t\geq 0$.

\begin{lemma}
For all $x,\lambda>0$, the function $(\mathcal{L}-\lambda)\widetilde{V}(x)=0$. 
\end{lemma}
{\bf Proof.} First note that we may take $f = \widetilde{V}$ in (\ref{itodecomp}) on account of the fact that $\widetilde{V}\in C^1(0,\infty)$. Since, for all $c>0$,
\[
e^{-\lambda t}\widetilde{V}(Z^c_t) -\widetilde{V}(c) - M^{\widetilde{V}}_t = \int_0^t e^{-\lambda s}(\mathcal{L}-\lambda)\widetilde{V}(Z^c_{s-}){\rm d}s, \,\,t\geq 0
\]
is a martingale, it must follow that, with probability one, $(\mathcal{L}-\lambda)\widetilde{V}(Z^c_t) = 0$ for Lebesgue almost every $t\geq 0$. It is a straightforward exercise however to deduce that $(\mathcal{L}- \lambda)\widetilde{V}(x)$ is a continuous function on $(0,\infty)$. Hence, since $Z^c_t$ is right continuous as a function of $t$, then so is $(\mathcal{L}-\lambda)\widetilde{V}(Z^c_t) $. It follows that, with probability one, $(\mathcal{L}-\lambda)\widetilde{V}(Z^c_t) = 0$ for all $t\geq 0$. In particular, for any $\theta\geq c>0$, $(\mathcal{L}-\lambda)\widetilde{V}(Z^c_{\tau^c_{(\theta, \infty)}}) =(\mathcal{L}-\lambda)\widetilde{V}(\theta) = 0$. This establishes the result.
\QED

Using the above lemma, we consider the behaviour of $(\mathcal{L}-\lambda)V^*(x)$ in the final step to proving the supermartingale property.
\begin{lemma}
For all $x>0$ we have $(\mathcal{L}-\lambda)V^*(x)\leq 0$ and hence it follows from (\ref{itodecomp}) that $\{e^{-\lambda t}V^*(Z^c_t): t\geq 0\}$ is a supermartingale.
\end{lemma}

{\bf Proof.} Note that for $x>b^*$ we have $V^*(x) = x$. It follows that for $x>b^*$
\[
(\mathcal{L}-\lambda)V^*(x) = (1- qx)+  \int_{(0,\infty)} \left\{V^*(e^{-\gamma y}x)-x\right\}m({\rm d}y)
\]
where we have used that $\lambda = q+\gamma\theta$. Note that for each $x>b^*$, for a given $\delta <x-b^*$ there exists an $\epsilon>0$ such that for all $z\in(x-\delta, x+\delta)$
\[
\int_{(0,\infty)} \left\{V^*(e^{-\gamma y}z)-z\right\}m({\rm d}y) = -z \int_{(0,\epsilon)} (1-e^{-\gamma y})m({\rm d}y) + \int_{[\epsilon,\infty)} \left\{V^*(e^{-\gamma y}z)-z\right\}m({\rm d}y).
\]
Using the above representation together with the fact that, by convexity, $V^{*\prime}(x)\leq 1$, a straightforward calculation shows that for all $x>b^*$
\[
\frac{\rm d}{{\rm d}x}(\mathcal{L}-\lambda)V^*(x) \leq-q - \int_{(0,\infty)}(1- e^{-\gamma y})m({\rm d}y)<0. 
\]
Note however that, thanks to the continuous and smooth pasting condition, 
$$(\mathcal{L}-\lambda)V^*(b^*) =(\mathcal{L}-\lambda)\widetilde{V}(b^*) = 0$$
 and hence it follows that 
$
(\mathcal{L}-\lambda)V^*(x) \leq 0
$
for all $x>b^*$. Since $(\mathcal{L}-\lambda)V^*(x)  = (\mathcal{L}-\lambda)\widetilde{V}(x)$ for $x\leq b^*$ we have that 
\begin{equation}
(\mathcal{L}-\lambda)V^*(x) \leq 0
\label{in hand}
\end{equation} 
for all $x>0$.

Now recalling that $V^*$ is a continuously differentiable function, using (\ref{itodecomp}) with $f = V^*$, we see, with the help of (\ref{in hand}), that for all $c>0$, $\{e^{-\lambda t}V^*(Z^c_t): t\geq 0\}$ is a supermartingale. Right continuity follows immediately by virtue of the fact that $V^*$ is continuous and that $\{Z^c_t : t\geq 0\}$ has right continuous paths.
 \QED

\begin{singlespace}

\end{singlespace}

\end{doublespace}
\end{document}